\magnification=\magstep1
\hsize=16.5 true cm 
\vsize=23.6 true cm
\font\bff=cmbx10 scaled \magstep1
\font\bfff=cmbx10 scaled \magstep2

\parindent0cm
%%%%%%%%%%%%%%%%%%%%%%%%%%%%%%%%
%%%%ABKšRZUNGEN im Textsatz %%%%
\def\cl{\centerline}           %
\def\bp{\bigskip}              %
\def\mp{\medskip}              %
\def\sp{\smallskip}            %
%%%%%%%%%%%%%%%%%%%%%%%%%%%%%%%%
%%%ZAHLENMENGEN%%%%%%%%%%%%%%%%%
           %
\def\Bbb#1{\hbox{\boldmas #1}} %
\def\R{\Bbb R}                 %
\def\N{\Bbb N}                 %
\def\Z{\Bbb Z}                 %
\def\bc{{\bf c}}
%%%%%%%%%%%%%%%%%%%%%%%%%%%%%%%%

%%%%%%%%%%%%%%%%%%%%%%%%%%%%%%%%
{\expandafter\edef\csname amssym.def\endcsname{%
       \catcode`\noexpand\@=\the\catcode`\@\space}

\catcode`\@=11

\def\undefine#1{\let#1\undefined}
\def\newsymbol#1#2#3#4#5{\let\next@\relax
 \ifnum#2=\@ne\let\next@\msafam@\else
 \ifnum#2=\tw@\let\next@\msbfam@\fi\fi
 \mathchardef#1="#3\next@#4#5}
\def\mathhexbox@#1#2#3{\relax
 \ifmmode\mathpalette{}{\m@th\mathchar"#1#2#3}%
 \else\leavevmode\hbox{$\m@th\mathchar"#1#2#3$}\fi}
\def\hexnumber@#1{\ifcase#1 0\or 1\or 2\or 3\or 4\or 5\or 6\or 7\or 8\or
 9\or A\or B\or C\or D\or E\or F\fi}

\font\tenmsa=msam10
\font\sevenmsa=msam7
\font\fivemsa=msam5
\newfam\msafam
\textfont\msafam=\tenmsa
\scriptfont\msafam=\sevenmsa
\scriptscriptfont\msafam=\fivemsa
\edef\msafam@{\hexnumber@\msafam}
\mathchardef\dabar@"0\msafam@39
\def\dashrightarrow{\mathrel{\dabar@\dabar@\mathchar"0\msafam@4B}}
\def\dashleftarrow{\mathrel{\mathchar"0\msafam@4C\dabar@\dabar@}}

\def\ulcorner{\delimiter"4\msafam@70\msafam@70 }
\def\urcorner{\delimiter"5\msafam@71\msafam@71 }
\def\llcorner{\delimiter"4\msafam@78\msafam@78 }
\def\lrcorner{\delimiter"5\msafam@79\msafam@79 }
\def\yen{{\mathhexbox@\msafam@55}}
\def\checkmark{{\mathhexbox@\msafam@58}}
\def\circledR{{\mathhexbox@\msafam@72}}
\def\maltese{{\mathhexbox@\msafam@7A}}

\font\tenmsb=msbm10
\font\sevenmsb=msbm7
\font\fivemsb=msbm5
\newfam\msbfam
\textfont\msbfam=\tenmsb
\scriptfont\msbfam=\sevenmsb
\scriptscriptfont\msbfam=\fivemsb
\edef\msbfam@{\hexnumber@\msbfam}
\def\Bbb#1{{\fam\msbfam\relax#1}}
\def\widehat#1{\setbox\z@\hbox{$\m@th#1$}%
 \ifdim\wd\z@>\tw@ em\mathaccent"0\msbfam@5B{#1}%
 \else\mathaccent"0362{#1}\fi}

\def\widetilde#1{\setbox\z@\hbox{$\m@th#1$}%
 \ifdim\wd\z@>\tw@ em\mathaccent"0\msbfam@5D{#1}%
 \else\mathaccent"0365{#1}\fi}
\font\teneufm=eufm10
\font\seveneufm=eufm7
\font\fiveeufm=eufm5
\newfam\eufmfam
\textfont\eufmfam=\teneufm
\scriptfont\eufmfam=\seveneufm
\scriptscriptfont\eufmfam=\fiveeufm

\newsymbol\risingdotseq 133A
\newsymbol\fallingdotseq 133B
\newsymbol\complement 107B
\newsymbol\nmid 232D
\newsymbol\rtimes 226F
\newsymbol\thicksim 2373

\font\eightmsb=msbm8   \font\sixmsb=msbm6   \font\fivemsb=msbm5
\font\eighteufm=eufm8  \font\sixeufm=eufm6  \font\fiveeufm=eufm5
\font\eightrm=cmr8     \font\sixrm=cmr6     \font\fiverm=cmr5
\font\eightbf=cmbx8    \font\sixbf=cmbx6    
      \font\eighti=cmmi8   \font\sixi=cmmi6
\font\ninesy=cmsy9     \font\eightsy=cmsy8  \font\sixsy=cmsy6
     \font\eightit=cmti8  
     \font\eightsl=cmsl8  
     \font\eighttt=cmtt8
 %SLANTED TYPEWRITER 10 POINT

\font\eightsmc=cmcsc8
\newskip\ttglue
\newfam\smcfam
\def\eightpoint{\def\rm{\fam0\eightrm}%
  \textfont0=\eightrm \scriptfont0=\sixrm \scriptscriptfont0=\fiverm
  \textfont1=\eighti \scriptfont1=\sixi \scriptscriptfont1=\fivei
  \textfont2=\eightsy \scriptfont2=\sixsy \scriptscriptfont2=\fivesy
  \textfont3=\tenex \scriptfont3=\tenex \scriptscriptfont3=\tenex
  \def\smc{\fam\smcfam\eightsmc}
  \textfont\smcfam=\eightsmc          
    %\scriptfont\smcfam=\sixsmc   \scriptscriptfont\smcfam=\fivesmc
\textfont\eufmfam=\eighteufm              \scriptfont\eufmfam=\sixeufm
     \scriptscriptfont\eufmfam=\fiveeufm
\textfont\msbfam=\eightmsb            \scriptfont\msbfam=\sixmsb
     \scriptscriptfont\msbfam=\fivemsb
\def\it{\fam\itfam\eightit}%
  \textfont\itfam=\eightit
  \def\sl{\fam\slfam\eightsl}%
  \textfont\slfam=\eightsl
  \def\bf{\fam\bffam\eightbf}%
  \textfont\bffam=\eightbf \scriptfont\bffam=\sixbf
   \scriptscriptfont\bffam=\fivebf
  \def\tt{\fam\ttfam\eighttt}%
  \textfont\ttfam=\eighttt
  \tt \ttglue=.5em plus.25em minus.15em
  \normalbaselineskip=9pt
  \def\MF{{\manual opqr}\-{\manual stuq}}%
  \let\big=\eightbig
  \setbox\strutbox=\hbox{\vrule height7pt depth2pt width\z@}%
  \normalbaselines\rm}
\def\eightbig#1{{\hbox{$\textfont0=\ninerm\textfont2=\ninesy
  \left#1\vbox to6.5pt{}\right.\n@space$}}}

\catcode`@=13 % reinstating the catcode of @

%%%%%%%%%%%%%%%%%%%%%%%%%%%%%%%%%%%%%%%

%%%%%%%%%%%%%%%%%%%%%%%%%%%%%%%%%%%%%%%
\cl{\bfff On compact subsets of the plane}
\mp
\centerline{\bff Gerald Kuba}
\bp
{\bff 1. Introduction}
\mp
Let $\,\bc\,$ denote the cardinality of the continuum. 
Let us call topological 
spaces $\;X_i\;(i\in I)\;$ {\it incomparable} if and only if 
$\,X_i\,$ is not homeomorphic to a subspace of $\,X_j\,$
whenever $\,i,j\in I\,$ and $\,i\not=j\,$.
An elementary but inconstructive proof shows that for any $\,n\in\N\,$
there exist $\,2^{\bc}\,$ incomparable subspaces of the Euclidean space
$\,\R^n\,$. (See [4] \S 35, V.~Theorem 1.)
Therefore and since it is clear that
there are precisely $\,\bc\,$ closed subsets of $\,\R^n\,$,
it is natural to ask whether 
there exist $\,\bc\,$ incomparable {\it compact} 
subspaces of $\,\R^n\,$.
This would be far from being true in the case $\,n=1\,$.
\mp
{\bf Proposition 1.} {\it Two compact subspaces 
of $\,\R\,$ are never incomparable.}
\mp
The situation depicted in Proposition 1 changes completely if 
the real line is replaced with the Euclidean space $\,\R^n\,$
for $\,n\geq 2\,$.
\mp
{\bf Theorem 1.} {\it There exist $\,\bc\,$ 
incomparable pathwise connected, compact subspaces of $\,\R^2\,$.}
\mp
Our goal is to prove Theorem 1 in a constructive and very natural way.
However, our proof is not short.
To put the theorem and its proof into perspective 
(and in view of the remark below) we show that
the following  counterpart of Theorem 1
can be verified in a very short and quite simple way.
\mp
{\bf Theorem 2.} {\it The plane $\,\R^2\,$ contains $\,\bc\,$ 
pathwise connected, nowhere dense, compact sets $\;X_i\;(i\in I)\;$
%with empty interiors 
such that for distinct $\,i,j\in I\,$
the spaces $\,X_i\,$ and $\,X_j\,$ are never homeomorphic,
whereas $\,X_i\,$ is homeomorphic to 
a subspace of $\,X_j\,$ whenever $\,i,j\in I\,$.}
\mp
{\it Remark.} Trivially, incomparable spaces are mutually non-homeomorphic. 
The well-known fact that 
$\,\R^n\,$ is homeomorphic to every open ball in $\,\R^n\,$
has two consequences. Firstly, the $\,\bc\,$ incomparable compact spaces in 
Theorem 1 must be {\it nowhere dense} subsets of $\,\R^2\,$.
%if in Theorem 2 no set $\,X_i\,$ were nowhere dense 
%then the last statement in Theorem~2 would hold trivially. 
Secondly, two {\it connected} subspaces of $\,\R\,$ are never 
incomparable.
\mp\sp
{\bff 2. Proof of Theorem 2}
\mp
For distinct points $\,a,b\in\R^2\,$ let $\,[a,b]\,$ denote the closed straight
line segment which connects $\,a\,$ with $\,b\,$
and put $\;]a,b]\,:=\,[a,b]\setminus\{a\}\,$.
For $\,n,k\in\N\,$ define 
\sp
\cl{$\;B(n)\,:=\;](2^{-n},0),(2^{-n},2^{-n})]\;$ and 
$\;T(n,k)\,:=\;](2^{-n},2^{-n}),(2^{-n}\!+\!3^{-nk},2^{-n+1})]\,$.} 
\sp
Let $\,{\cal M}\,$ be the family of all infinite sets of integers $\,k\geq 4\,$.
Naturally, the cardinality of 
$\,{\cal M}\,$ is $\,\bc\,$.
For $\,M\in{\cal M}\,$ define 
a pathwise connected subset $\,A[M]\,$ of $\,\R^2\,$ via 
\sp
\cl{$ A[M]\;:=\;([0,1]\!\times\!\{0\})\;\cup\!
\bigcup\limits_{n\in M}\!\big(\,B(n)\cup
\bigcup\limits_{k=1}^{n-1}T(n,k)\,\big)\,$.}
\sp
Obviously, $\,A[M]\,$ is compact and nowhere dense.
One may picture $\,A[M]\,$ as a {\it tree} where 
$\;[0,1]\!\times\!\{0\}\;$ is the {\it trunk} and 
$\;B(n)\;(n\in M)\;$ are the {\it branches} 
and $\,n-1\,$ {\it twigs} $\,T(n,k)\,$
are attached to the branch $\,B(n)\,$ for every $\,n\in{\cal M}\,$.
The $\,\bc\,$ spaces $\,A[M]\,$ are mutually non-homeomorphic 
because each set $\,M\in{\cal M}\,$ is completely determined 
by the topology of $\,A[M]\,$ via 
\sp
\cl{$\;M\,=\,\{\,p(x)\;|\;x\in A[M]\,\}\setminus\{1,2,3\}\;$}
\sp
where $\,p(x)\,$ is the number of the path-components 
of the space $\,A[M]\setminus\{x\}\,$.
\eject
\sp
Finally, if $\,M_1,M_2\in{\cal M}\,$ then
some subspace of $\,A[M_2]\,$ is homeomorphic to $\,A[M_1]\,$.
Indeed, since $\,M_2\,$ is an {\it infinite} subset of $\,\N\,$,
such a subspace can easily be obtained from the tree $\,A[M_2]\,$ 
by cutting off twigs or by eliminating branches 
together with all adjacent twigs.
\bp
{\bff 3. Proof of Theorem 1} 
\mp
As usual, $\,|S|\,$ denotes the cardinal number
of the set $\,S\,$. 
We will prove Theorem 1 in three main steps.
Firstly we construct $\,\bc\,$ 
incomparable connected, closed subspaces $\,Y\,$ of $\,\R^2\,$
which are neither compact nor pathwise connected.
Secondly, by adding straight lines we modify these spaces 
a bit such that incomparability is saved.
Thirdly, by adding an appropriate curve 
around the modified spaces we obtain
$\,\bc\,$ incomparable pathwise connected, compact subspaces $\,X\,$ 
of $\,\R^2\,$. 
The basic idea to carry out the first step is to consider a double sequence
\sp
\cl{$\dots,\,a_{-3},\,a_{-2},\,a_{-1},\,
a_{0},\,a_{1},\,a_{2},\,a_{3},\,\dots$}
\sp
of digits $\,a_k\in\{0,1\}\,$ and then construct 
a connected, closed subspace $\,Y\,$ of $\,\R^2\,$
such that the topology of $\,Y\,$
completely determines the double sequence 
up to shifts and reflections.
\mp
Let $\,\pi\,$ denote the projection 
$\;(x,y)\mapsto x\;$ from $\,\R^2\,$ onto $\,\R\,$.
Note that for any compact and connected set $\,A\subset\R^2\,$
the set $\,\pi(A)\,$ is a compact interval.
With the notation $\,[\cdot,\cdot]\,$ for compact straight lines in 
the plane as in Section 2,
for every $\,k\in\Z\,$ put 
\mp
\qquad $U_k\;:=\;[(k\!-\!{1\over 3},-1),(k\!-\!{1\over 3},1)]\,\cup\,
[(k\!-\!{1\over 3},-1),(k\!+\!{1\over 3},-1)]\,\cup\,
[(k\!+\!{1\over 3},-1),(k\!+\!{1\over 3},1)]$\quad and
\sp
\qquad $H_k\;:=\;[(k\!-\!{1\over 3},-1),(k\!-\!{1\over 3},1)]\,\cup\,
[(k\!-\!{1\over 3},0),(k\!+\!{1\over 3},0)]\,\cup\,
[(k\!+\!{1\over 3},-1),(k\!+\!{1\over 3},1)]\;$.
\mp
So $\,U_k\,$ is a curve which looks like the figure $\,\bigsqcup\,$
(and hence $\,U_k\,$ is homeomorphic with the unit interval $\,[0,1]\,$)
and $\,H_k\,$ is a point set                  
which looks like the figure $\;|\!\!\!-\!\!\!|\,$.
Obviously, 
$\,U_k\,$ and $\,H_k\,$ are pathwise
connected, compact subsets of $\;[k-{1\over 3},k+{1\over 3}]\times[-1,1]\;$
and $\;\pi(U_k)=\pi(H_k)=[k-{1\over 3},k+{1\over 3}]\,$.
\mp
Next, for every $\,k\in\Z\,$ define 
\mp
\cl{$\;S_k\;:=\;
\big\{\,\big(x,\,\sin\,((x-k-{1\over 3})^{-1}(x-k-{2\over 3})^{-1})\big)
\;\,\big|\,\;k+{1\over 3}\,<\,x\,<\,k+{2\over 3}\,\big\}\,.$}
\mp
Thus $\,S_k\,$ is a pathwise connected subset of 
$\;]k+{1\over 3},k+{2\over 3}[\,\times[-1,1]\;$
homeomorphic to $\,\R\,$ and $\;\pi(S_k)=\,]k+{1\over 3},k+{2\over 3}[\,$.
Naturally, the closure $\,\overline{S_k}\,$ of $\,S_k\,$ in $\,\R^2\,$
is connected and compact. More precisely, 
$\;\overline{S_k}\,=\,I_1\cup S_k\cup I_2\;$
with $\;I_n\,:=\,
[(k\!+\!{n\over 3},-1),(k\!+\!{n\over 3},1)]\;$ for $\,n\in\{1,2\}\,$.
Of course,  
$\,I_1\,$ and $\,S_k\,$ and $\,I_2\,$
are the path-components of the connected, compact space $\,\overline{S_k}\,$.
Obviously,
\sp
\cl{$\;\bigcup\limits_{k=-\infty}^\infty\pi(U_k)\,\;=\, 
\bigcup\limits_{k=-\infty}^\infty\pi(H_k)\,\;=\,\; 
\R\,\setminus\bigcup\limits_{k=-\infty}^\infty\pi(S_k)\;.$}
\mp
Let $\,{\cal G}_0\,$ denote the family of all subspaces of $\,\R^2\,$
which are homeomorphic to the unit interval $\,[0,1]\,$. 
Let $\,{\cal G}_1\,$ denote the family of all subspaces of $\,\R^2\,$
homeomorphic to the space $\,H_0\,$.
Let $\,{\cal G}_2\,$ denote the family of all subspaces of $\,\R^2\,$
homeomorphic to $\,\R\,$.
By definition, $\;U_k\in{\cal G}_0\;$
and $\;H_k\in{\cal G}_1\;$ and $\;S_k\in{\cal G}_2\;$
for every $\,k\in\Z\,$.
Clearly, if $\,X_n\in{\cal G}_n\,$ for $\,n\in\{0,1,2\}\,$ 
then $\,X_0\,$ and $\,X_1\,$ and $\,X_2\,$ are mutually  
non-homeomorphic. (This is elementary since $\,[0,1]\,$ 
is compact and has precisely two 
noncut points and $\,H_0\,$ is compact and has precisely four noncut points
and the real line $\,\R\,$ is not compact.)
\mp
Let $\,\Omega\,$ denote the family of all mappings $\,g\,$
from $\,\Z\,$ to $\,\{0,1\}\,$. For every $\,g\in\Omega\,$ define 
a subset $\,Y[g]\,$ of $\;\R\times[-1,1]\;$ 
with $\,\pi(Y[g])=\R\,$ via
\sp
\cl{$Y[g]\;\;:=\;\;\bigcup\limits_{k=-\infty}^{\infty}(S_k\cup A_k(g))$}
\sp
where $\,A_k(g)=U_k\,$ when $\,g(k)=0\,$ and   
$\,A_k(g)=H_k\,$ when $\,g(k)=1\,$.
Naturally, the set 
$\;\bigcup\limits_{k=-n}^{n}(A_k(g)\cup S_k\cup A_{k+1}(g))\;$
is connected and compact for every integer $\,n\geq 0\,$.
Therefore, $\,Y[g]\,$ is a connected and closed subset of $\,\R^2\,$. 
\mp
For $\,g\in\Omega\,$ let $\,{\cal P}[g]\,$ 
denote the family of all path-components
of $\,Y[g]\,$. Obviously, if 
$\;{\cal P}_n[g]\,:=\,{\cal P}[g]\cap{\cal G}_n\;$ for $\,n\in\{0,1,2\}\,$
then 
\sp
\cl{$\;{\cal P}[g]\,=\,{\cal P}_0[g]\cup{\cal P}_1[g]\cup{\cal P}_2[g]\;$}
\sp
and the subfamilies $\;{\cal P}_0[g],{\cal P}_1[g],{\cal P}_2[g]\;$
of $\,{\cal P}[g]\,$ are mutually disjoint.
By definition, $\;{\cal P}_0[g]\,=\,\{\,U_k\;|\;k\in g^{-1}(\{0\})\,\}\;$ 
and $\;{\cal P}_1[g]\,=\,\{\,H_k\;|\;k\in g^{-1}(\{1\})\,\}\;$ 
and $\;{\cal P}_2[g]\,=\,\{\,S_k\;|\;k\in\Z\,\}\,$. 
Furthermore, the following is clearly true. 
\mp
(3.1)$\;$ {\it If $\;P_n\in{\cal P}_n[g]\;$ for $\,n\in\{0,1,2\}\;$ 
and $\,i,j\in\{0,1,2\}\,$ 
then $\,P_i\,$ is homeomorphic with $\,P_j\,$
if and only if $\,i=j\,$.}
\mp
The family $\,{\cal P}[g]\,$
is naturally ordered
by defining $\,P_1\prec P_2\,$ for $\,P_1,P_2\in{\cal P}[g]\,$
if and only if $\,x_1<x_2\,$ for some $\,(x_1,y_1)\in P_1\,$ 
and some $\,(x_2,y_2)\in P_2\,$.
In particular, the subfamily 
$\;{\cal P}^*[g]\,:=\,{\cal P}_0[g]\cup{\cal P}_1[g]\;$
is linearly ordered and for $\,A, B\in{\cal P}^*[g]\,$
we have $\,A\prec B\,$ if and only if the unique integer
in $\,\pi(A)\,$ is smaller than the unique integer
in $\,\pi(B)\,$.
Of course, the linearly ordered set $\,({\cal P}[g],\prec)\,$ 
and the linearly ordered set $\,({\cal P}^*[g],\prec)\,$  
are both isomorphic with the naturally ordered set $\,\Z\,$. 
\mp
Obviously, if $\,P_1,P_2\in{\cal P}[g]\,$ and $\,P_1\prec P_2\,$
then $\,P_1\cup P_2\,$ is connected if and only if  
there does not exist a $\,P\in{\cal P}[g]\,$
with $\,P_1\prec P\prec P_2\,$.
This has the following important consequence. 
If $\,g_1,g_2\in\Omega\,$ and 
$\,\varphi\,$ is a homeomorphism from $\,Y[g_1]\,$ onto $\,Y[g_2]\,$
then $\,P\mapsto \varphi(P)\,$ is a monotonic bijection 
from the linearly ordered set $\,{\cal P}[g_1]\,$ onto 
the linearly ordered set $\,{\cal P}[g_2]\,$. 
In view of (3.1) we have $\;\varphi(P)\in{\cal P}_n[g_2]\;$
whenever $\;P\in{\cal P}_n[g_1]\;$ and $\,n\in\{0,1,2\}\,$.
Therefore, $\,\varphi\,$ induces 
a monotonic bijection $\,\psi\,$
from the linearly ordered set $\,{\cal P}^*[g_1]\,$ onto 
the linearly ordered set $\,{\cal P}^*[g_2]\,$ in view of (3.1).
Suppose that $\,\psi\,$ is increasing.
Since $\;\varphi(P)\in{\cal P}_n[g_2]\;$
whenever $\;P\in{\cal P}_n[g_1]\;$ and $\,n\in\{0,1\}\,$,
we conclude that the following statement is true.
\sp
(3.2)$\;${\it There exists an integer $\,m\,$
such that $\;g_1(m+k)=g_2(k)\;$ for every $\,k\in\Z\,$.}
\sp
On the other hand, if $\,\psi\,$ is decreasing 
then the following is true.
\sp
(3.3)$\;${\it There exists an integer $\,m\,$
such that $\;g_1(m-k)=g_2(k)\;$ for every $\,k\in\Z\,$.}
\sp
So if $\,g_1,g_2\,$ are chosen in $\,\Omega\,$
such that neither (3.2) nor (3.3) holds then we can be sure that the spaces
$\,Y[g_1]\,$ and $\,Y[g_2]\,$ are not homeomorphic.
Therefore we can track down $\,\bc\,$ 
mutually non-homeomorphic spaces $\,Y[g]\,$ with $\,g\in\Omega\,$ 
as follows.
\mp
Let $\,\Omega^*\,$ be the family of all double sequences $\,g\in\Omega\,$
such that $\,g(0)=0\,$ and $\,g(k)=1\,$ for every negative integer $\,k\,$
and $\,g(k)=0\,$ for infinitely many positive integers $\,k\,$.
Naturally, $\,|\Omega^*|=\bc\,$. It is evident that 
for distinct $\,g_1,g_2\in\Omega^*\,$ neither (3.2) nor (3.3) is true
and hence the $\,\bc\,$ spaces $\;Y[g]\;(g\in\Omega^*)\;$
are mutually non-homeomorphic. We claim that these spaces are 
also incomparable. To prove this, now it is enough to verify that for 
$\,g,h\in\Omega\,$ the space $\,Y[g]\,$ cannot be homeomorphic with 
a {\it proper} subspace of $\,Y[h]\,$. 
\mp
Assume indirectly that $\,f\,$ is a homeomorphism from $\,Y[g]\,$
onto a subspace of $\,Y[h]\,$ such that $\,f(Y[g])\not=Y[h]\,$.
Since $\,Y[g]\,$ is connected, $\,f(Y[g])\,$ is connected.
Furthermore, for every $\,P\in{\cal P}[g]\,$ 
the set $\,f(P)\,$ is a path-component in the space $\,f(Y[g])\,$.
Therefore, for every $\,P\in{\cal P}[g]\,$ 
there is a unique $\,Q\in{\cal P}[h]\,$ such that
$\,f(P)\subset Q\,$. Of course, $\,f\,$ induces a 
linear ordering of the family $\,{\cal C}\,$ of all path-components 
of the space $\,f(Y[g])\,$. So $\,{\cal C}\,$ is linearly ordered
as $\,\Z\,$. Consequently, $\,\pi(f(Y[g]))=\R\,$.
Furthermore, since
every $\,C\in{\cal C}\,$ has 
an immediate predecessor and an immediate successor in $\,{\cal C}\,$,
we can be sure that $\,f(P)=Q\,$ and 
$\,Q\in{\cal P}_2[h]\,$ if $\,P\in{\cal P}_2[g]\,$
because if $\,f(P)\not=Q\,$ then the space $\,f(Y[g])\,$ 
would not be connected. Therefore we can be sure that the following
statement is true for $\,n=2\,$.
\sp
(3.4)$\;$ {\it $\,P\mapsto f(P)\,$ defines an injective mapping 
from $\,{\cal P}_n[g]\,$ into $\,{\cal P}_n[h]\,$.}
\sp
Moreover, for $\,n=0\,$ the mapping in (3.4) must be 
bijective since $\,\pi(f(Y[g]))=\R\,$.
In order to verify (3.4) also for $\,n=1\,$ and $\,n=2\,$,
we firstly point out that it is impossible 
that $\,f(P)\subset Q\,$ for $\,P\in{\cal P}_1[g]\,$
and $\,Q\in{\cal P}_0[h]\,$. Because it is clear that 
a space in $\,{\cal G}_1\,$ cannot be homeomorphic to 
a subspace of a space in $\,{\cal G}_0\,$.
The converse is not true because 
it is obvious that 
$\,U_0\,$ is homeomorphic to a subspace of $\,H_0\,$.
Nevertheless, $\,f\,$ cannot map any set $\,P\in{\cal P}_0[g]\,$
onto a (necessarily connected) 
subset of a set from the family $\,{\cal P}_1[h]\,$.
Because otherwise we clearly could track down an infinite 
subset $\,\Delta\,$  of one set $\,S_k\in{\cal G}_2\,$ with 
$\,S_k\subset f(Y[g])\,$
(where $\,S_k\,$ is either an immediate predecessor or an immediate successor 
of $\,f(P)\in{\cal C}\,$)
such that $\,\Delta\,$ is a discrete subset of $\,f(Y[g])\,$.
This is impossible because the closure of $\,f^{-1}(S_k)\,$
in $\,Y[g]\,$ is compact and hence 
the infinite set $\,f^{-1}(\Delta)\,$ cannot be discrete in 
the space $\,Y[g]\,$. For the same reason 
it is impossible that for $\,n\in\{1,2\}\,$ 
and $\,P\in{\cal P}_n[g]\,$ the path-component 
$\,f(P)\,$ in the space $\,f(Y[g])\,$ is a 
proper subset of some set in the family $\,{\cal P}_n[h]\,$.
Thus we conclude that
$\,(3.4)\,$ is true for every $\,n\in\{0,1,2\}\,$.
\mp
Consequently, in view of $\,\pi(f(Y[g]))=\R\,$ 
we arrive at $\,f(Y[g])=Y[h]\,$ which contradicts 
our indirect assumption about the mapping $\,f\,$.
Thus we can be sure that the $\,\bc\,$ spaces $\;Y[g]\;(g\in\Omega^*)\;$
are incomparable. This concludes the first main step in the proof 
of Theorem 1.
\mp
In the second main step we turn the $\,\bc\,$ connected spaces $\,Y[g]\,$
into $\,\bc\,$ connected but not pathwise connected
subspaces $\,X[g]\,$ of $\,\R^2\,$
such that for some point $\,z\in\R^2\,$ 
the set $\,X[g]\cup\{z\}\,$ is pathwise connected. 
\mp
Let $\,\phi\,$ denote the mapping 
$\;(x,y)\mapsto ({2\over \pi}\arctan x,\,y)\;$ which, of course, 
is a homeomorphism from $\,\R^2\,$
onto $\;]{-1,1}[\times\R\,$. 
Therefore, for every $\,g\in\Omega\,$
the space $\,Y[g]\,$ is homeomorphic to the space
$\,\phi(Y[g])\,$. In particular, the $\,\bc\,$ 
connected spaces $\;\phi(Y[g])\;(g\in\Omega^*)\;$ are incomparable. 
But as subsets of $\,\R^2\,$ the set $\,Y[g]\,$ and the set
$\,\phi(Y[g])\,$ are of different kind. 
While $\,Y[g]\,$ is closed and unbounded, 
$\,\phi(Y[g])\,$ is bounded but not closed.
By definition, $\;\phi(Y[g])\,\subset\;]{-1,1}[\times[-1,1]\;$ 
for every $\,g\in\Omega\,$. 
\mp
Since $\;U_k\cap H_k\cap(\R\times\{1\})\,=\,
\{(k\!-\!{1\over 3},1),\,(k\!+\!{1\over 3},1)\}\;$
for every $\,k\in\Z\,$, we can fix a set 
$\;\Gamma\,\subset\,\R\times\{1\}\;$ such that for arbitrary $\,g\in\Omega^*\,$
every path-component of $\,\phi(Y[g])\,$ contains precisely one point 
$\,(x,1)\in \Gamma\,$ and every point $\,(x,1)\in \Gamma\,$ lies in 
some path-component of $\,\phi(Y[g])\,$.
For every $\,a\in \Gamma\,$ consider the compact straight line 
segment $\;[a,(1,2)]\;$
which connects the point $\,a\,$ with the 
point $\,(1,2)\,$ and put $\;L_a\,:=\,[a,(1,2)]\setminus\{(1,2)\}\,$.
For every $\,g\in\Omega^*\,$ put 
\sp
\cl{$\;X[g]\;:=\;\phi(Y[g])\cup\bigcup\limits_{a\in \Gamma}L_a\,$.}
\sp
Certainly, $\,X[g]\cup\{(1,2)\}\,$ is a 
pathwise connected subspace of $\,\R^2\,$
while $\,X[g]\,$ is a connected space with infinitely many 
path-components since $\,\phi(Y[g])\subset X[g]\,$.
\mp
For $\,g\in\Omega^*\,$ let
$\,D\,$ be a path-component in the connected space $\,X[g]\,$.
Then there is a unique path-component $\,P_D\,$ in $\,Y[g]\,$
such that $\;D\,=\,\phi(P_D)\cup L_a\;$
for $\,a=a_D\,$ defined by $\,\{a\}=\phi(P_D)\cap\Gamma\,$. If $\,P_D\,$
is compact then either $\,P_D=U_k\,$ or $\,P_D=H_k\,$ for some $\,k\in\Z\,$
and hence it is plain that $\,D\,$ is homeomorphic 
with $\,P_D\setminus\{\phi^{-1}(a_D)\}\,$. 
Clearly, $\,D\,$ has precisely three noncut points
if $\,P_D=H_k\,$ and precisely one noncut point if $\,P_D=U_k\,$.
For the moment let 
us call a point $\,z\,$ in a pathwise connected space $\,Z\,$ 
a {\it triple point}
when $\,Z\setminus\{z\}\,$ has precisely {\it three} path-components. 
If $\,P_D\,$ is not compact then 
$\;P_D=S_k\;$ for some $\,k\in\Z\,$ and hence the point $\,a_D\,$ 
is the unique triple point of $\,D\,$.
On the other hand, if $\,P_D\,$ is compact then $\,D\,$
has either precisely two triple points or no triple point.
\mp
In view of these considerations we realize that 
by adding straight line segments $\,L_a\,$ to the path-components
of $\,\phi(Y[g])\,$ we do not create an essentially new situation
and that, by exactly the same arguments which prove the $\,\bc\,$
spaces $\;Y[g]\;$ to be incomparable,
the connected spaces $\;X[g]\;(g\in\Omega^*)\;$ must be
incomparable as well. (By the way, the following considerations 
show that the $\,\bc\,$ pathwise connected spaces $\;X[g]\cup\{(1,2)\}\;\;$
are incomparable too.)
\mp
Now, in the final step of the proof we bring compactness into play.
Consider the trapezoidal curve
\sp
\cl{$Q\;:=\;[(-1,-1),(1,-2)]\cup[(1,-2),(1,2)]
\cup[(1,2),(-1,1)]\cup[(-1,1),(-1,-1)]$}
\sp
and for every $\,g\in\Omega^*\,$ put 
\sp
\cl{$C[g]\;:=\;X[g]\cup Q\;=\;
\phi(Y[g])\cup Q\cup\bigcup\limits_{a\in \Gamma}[a,(1,2)]\;$.}
\mp
It is evident that $\,C[g]\,$ is compact and pathwise connected
for every $\,g\in\Omega^*\,$. We conclude the proof by 
verifying that the spaces $\;C[g]\;(g\in\Omega^*)\;$
are incomparable. 
\mp
Let $\,g_1,g_2\in\Omega^*\,$ and assume that $\,f\,$ is a homeomorphism 
from $\,C[g_1]\,$ onto a subspace of $\,C[g_2]\,$.
For the moment let us call a point $\,b\,$ in a 
pathwise connected space $\,X\,$ 
a {\it super path-cut point} 
when $\,X\setminus\{b\}\,$ has infinitely many path-components. 
Obviously, the vertex $\,(1,2)\,$ of the trapezoid $\,Q\,$ is 
a super path-cut point in $\,C[g_1]\,$ and it is the only one. 
Thus the pathwise connected space $\,f(C[g_1])\,$ has precisely one 
super path-cut point $\,b\,$ and it is evident that 
there is no other possibility than $\,b=(1,2)\,$.
Hence, $\,(1,2)\,$ is a fixed point of the mapping $\,f\,$.
In particular, $\;(1,2)\in f(C[g_1])\,$.
Trivially, the closure of the set $\,Q\setminus\{(1,2)\}\,$ in 
the space $\,C[g_i]\,$ is $\,Q\,$ for $\,i=1\,$ and $\,i=2\,$.
Obviously, $\,Q\setminus\{(1,2)\}\,$ is the unique 
path-component of the connected space 
$\;C[g_1]\setminus\{(1,2)\}\;$ whose closure in $\,C[g_1]\,$
is a connected space without cut points.
Moreover, if $\,A\subset C[g_i]\,$ is pathwise connected
and $\,(1,2)\not\in A\,$ 
and the closure of $\,A\,$
in the space $\,C[g_i]\,$ 
is a connected space without cut points
then there is no other possibility than 
$\;A\,=\,Q\setminus\{(1,2)\}\,$.  
Therefore, $\,f(Q)=Q\,$ and hence via $\;C[g_i]\setminus Q=X[g_i]\;$
we derive $\;f(X[g_1])\subset X[g_2]\;$ 
and hence $\,g_1=g_2\,$ since $\,X[g_1],X[g_2]\,$ are incomparable 
spaces.
\bp
{\it Remark.} The $\,\bc\,$ subspaces $\;Q\cup\phi(Y[g])\;$ of 
$\,\R^2\,$ with $\,g\in\Omega^*\,$ 
are obviously {\it connected and compact}. They are also {\it incomparable}
because $\,Q\,$ is clearly the unique infinite pathwise connected 
subspace of $\;Q\cup\phi(Y[g])\;$ without cut points 
for every $\,g\in\Omega^*\,$
and the spaces $\;\phi(Y[g])\;(g\in\Omega^*)\;$ are incomparable. 
So the proof of Theorem 1 can be abridged if the property 
{\it pathwise connected} is replaced by the weaker property
{\it connected}. Concerning {\it pathwise connectedness} 
without compactness,
it is quite easy to find a family $\,{\cal F}\,$
of incomparable, pathwise connected subspaces of $\,\R^2\,$
with $\,|{\cal F}|=2^\bc\,$.
Indeed, apply [4] \S 35, V.~Theorem 1 in order to define 
a family $\,{\cal G}\,$ of $\,2^\bc\,$ incomparable 
(and hence totally disconnected) 
subspaces of $\,\R\,$. Then, with 
$\;X_G\,:=\,(\R\times\{0\})\cup(G\times[0,1])\;$
for each $\,G\in{\cal G}\,$,
the family $\;{\cal F}\,:=\,\{\,X_G\;|\;G\in{\cal G}\,\}\;$
does the job.  
(Because for $\,G\in{\cal G}\,$ and any pathwise connected subspace 
$\,Z\,$ of $\,X_G\,$ all triple points in $\,Z\,$ must lie
in $\,G\times\{0\}\,$.)
\bp
{\bff 4. Proof of Proposition 1} 
\mp
Assume that $\,\emptyset\not=X\subset\R\,$ 
and $\,X\,$ is not totally disconnected.
Then $\;]a,b[\,\subset X\;$ for some $\,a<b\,$ and,
since $\,\R\,$ and $\,]a,b[\,$ are homeomorphic,
every subspace of $\,\R\,$ is homeomorphic to a subspace of $\,X\,$. 
Therefore, Proposition 1 is an immediate consequence 
of the following noteworthy proposition.
\mp
{\bf Proposition 2.} {\it If $\,{\cal F}\,$ is a family of incomparable
totally disconnected, compact, metrizable
spaces then $\,|{\cal F}|\leq 1\,$.}
\mp
{\it Proof.} We apply well-known basic facts from topology (see [4]).
Every {\it totally disconnected} compact Hausdorff space is 
{\it zero-dimensional.} Every compact, metrizable space 
is separable and completely 
metrizable.
If $\,X\,$ is a
zero-dimensional, separable, completely metrizable space
and $\,{\Bbb D}\,$ denotes the Cantor set
then $\,X\,$ is homeomorphic to some subspace of $\,{\Bbb D}\,$.
Furthermore, some subspace of $\,X\,$ is homeomorphic to $\,{\Bbb D}\,$
provided that $\,X\,$ is uncountable.
Consequently, $\,|{\cal F}|=1\,$ if some $\,X\in{\cal F}\,$ is uncountable. 
\sp
It is well-known (see [1] Corollary 1095)
that the topology of any compact, countable Hausdorff space $\,H\,$
is the order topology of some well-ordering of $\,H\,$.
Naturally if $\,(X_i,\preceq_i)\,$
are well-ordered sets for $\,i\in\{1,2\}\,$ then 
$\,X_1\,$ and $\,X_2\,$ are isomorphic or
$\,X_i\,$ is isomorphic with $\;\{\,x\in X_j\;|\;x\preceq_j y\,\}\;$ 
for some $\,y\in X_j\,$ where either 
$\,(i,j)=(1,2)\,$ or $\,(i,j)=(2,1)\,$.
Consequently, two compact, countable Hausdorff spaces are never 
incomparable. In particular, $\,|{\cal F}|\leq 1\,$
if all spaces in $\,{\cal F}\,$ are countable, {\it q.e.d.}
\bp
{\it Remark.} Consider the Euclidean space $\,\R^n\,$ for any 
dimension $\,n\,$.
The previous proof explains that if 
$\;X_i\;(i\in I)\;$ are uncountable, totally disconnected, compact
subspaces of $\,\R^n\,$ then $\,X_i\,$ is homeomorphic to 
a subspace of $\,X_j\,$ whenever $\,i,j\in I\,$.
Therefore, in comparison with Theorem 2 it is worth mentioning
that there exist $\,\bc\,$ mutually non-homeomorphic 
totally disconnected, compact subspaces $\,X\,$ of $\,\R^n\,$
with $\,|X|=\bc\,$. (See [2] Proposition 2.)
\bp\mp
{\bff 5. A stronger version of Theorem 1}
\mp
The $\,\bc\,$ spaces %constructed 
in the proof of Theorem 2 are obviously {\it locally 
pathwise connected}.
On the other hand, due to the building blocks $S_k$,
none of the $\,\bc\,$ spaces proving Theorem 1 is 
{\it locally connected}. There arises a natural question 
which is answered by the following theorem.
\bp 
{\bf Theorem 3.} {\it There exist $\,\bc\,$ 
incomparable pathwise connected and locally connected 
compact subspaces of $\,\R^2\,$.}
\mp
A proof of Theorem 3 can be found in Section 6 of our paper [3].
(In Theorem 3 the word {\it pathwise} could be omitted 
since every connected, locally connected, 
compact and metrizable space is pathwise connected.)
\bp\bp
{\bff References}
\medskip
[1] Dasgupta, A.: Set Theory,  Birkh\"auser 2014.
\smallskip
[2] Kuba, G.: {\it Counting ultrametric spaces.}
Colloq.~Math. {\bf 152} (2018), 217-234.
\sp
[3] Kuba, G.: {\it Proper classes of non-embeddable continua.}
{\tt arXiv:2609.17173 (2026)}
\sp
[4] Kuratowski, K.: Topology, vol.~1. Academic Press 1966. 
\bp\bp
{\sl Author's address:} Institute of Mathematics. 

University of Natural Resources and Life Sciences, Vienna, Austria. 
\sp
{\sl E-mail:} {\tt gerald.kuba(at)boku.ac.at}
\end